\documentclass[a4paper,12pt]{amsart}

\usepackage{float}
\usepackage{euscript,eufrak,verbatim}
\usepackage{graphicx}
\usepackage[usenames]{color}
\usepackage[colorlinks,linkcolor=red,anchorcolor=blue,citecolor=blue]{hyperref}
\usepackage{amsmath}
\usepackage{amsthm}
\usepackage[all]{xy}
\usepackage{graphicx} 
\usepackage{amssymb}

\theoremstyle{plain}
\newtheorem{main theorem}{Main Theorem}
\newtheorem{theorem}{Theorem}[section]
\newtheorem{lemma}[theorem]{Lemma}

\newtheorem{corollary}[theorem]{Corollary}
\newtheorem{proposition}[theorem]{Proposition}
\newtheorem{claim}[theorem]{Claim}

\theoremstyle{definition}

\newtheorem{problem}[theorem]{Problem}

\makeatletter 
\@addtoreset{equation}{section}
\numberwithin{equation}{section}

\newcommand{\norm}[1]{\left\lVert#1\right\rVert}

\newcommand{\diam}{\mathrm{Diam}}

\newcommand{\mdim}{\mathrm{mdim}}

\newcommand{\widim}{\mathrm{Widim}}
\newcommand{\htop}{h_{\mathrm{top}}}

\newcommand{\lmdimm}{\underline{\mathrm{mdim}}_{\mathrm{M}}}

\title[Waist inequality, entropy and mean dimension: II]{Application of waist inequality to entropy and mean dimension: II}

\author{Masaki Tsukamoto}

\address
{Department of Mathematics, Kyoto University, Kitashirakawa Oiwake-cho, Sakyo-ku, Kyoto 606-8502, Japan}

\email{tsukamoto@math.kyoto-u.ac.jp}

\begin{document}



\thanks{M.T. was supported by JSPS KAKENHI JP21K03227.}

\maketitle

\begin{abstract}
Let $X$ be a full-shift on the alphabet $[0, 1]^a$ and let $(Y, S)$ be an arbitrary dynamical system.
We prove that any equivariant continuous map from $X$ to $Y$ has conditional metric mean dimension not 
less than $a-\mdim(Y, S)$.
This solves a problem posed in \cite[Problem 1.6]{Shi--Tsukamoto}.
\end{abstract}

\section{Main result}  \label{section: main result}

This paper is a continuation of \cite{Shi--Tsukamoto}.
Our purpose is to solve a problem left open in \cite[Problem 1.6]{Shi--Tsukamoto}.

In \cite{Shi--Tsukamoto} we developed applications of waist inequality to topological dynamics.
Waist inequality is a fundamental inequality in geometry and topology
discovered by Gromov \cite{Gromov83, Gromov03}.
It states that any continuous map from the sphere $S^n$ to $\mathbb{R}^m$ $(n\geq m)$ must have a “large” fiber.
Waist inequality and its variants have appeared in many interesting subjects 
\cite{Giannopoulos--Milman--Tsolomitis, Vershynin, Gromov09, Gromov10, Gromov--Guth, Klartag17, Klartag18}.

The paper \cite{Shi--Tsukamoto} discovered applications of waist inequality to entropy and mean dimension of dynamical systems.
We quickly prepare basic notions and review a result of \cite{Shi--Tsukamoto}.

Let $(X, d)$ be a compact metric space.
For a positive number $\varepsilon$ we define $\#(X, d, \varepsilon)$ as the minimum number $n$ such that 
there exists an open cover $X= U_1\cup U_2\cup \dots \cup U_n$ with $\diam U_i < \varepsilon$ for all $i$.
A continuous map $f\colon X\to Y$ to a topological space $Y$ is called an \textbf{$\varepsilon$-embedding}
if $\diam f^{-1}(y) < \varepsilon$ for all $y\in Y$.
We define $\widim_\varepsilon(X, d)$ as the minimum number $n$ such that there exists a finite simplicial complex $P$ of dimension $n$ and 
an $\varepsilon$-embedding $f\colon X\to P$.
(Here "finite" means that the number of vertices of $P$ is finite.)

A pair $(X, T)$ is called a \textbf{dynamical system} if $X$ is a compact metrizable space and $T\colon X\to X$ is a 
homeomorphism.
Given a dynamical system $(X, T)$, we will define its topological entropy and mean dimension.
Let $d$ be a metric on $X$. For each natural number $N$ we define a new metric $d_N$ on $X$ by 
\[ d_N(x, y) = \max_{0\leq n <N} d(T^n x, T^n y). \]
We define the \textbf{topological entropy} $\htop(X, T)$ and \textbf{mean dimension} $\mdim(X, T)$ by 
\begin{align*}
   \htop(X, T) &= \lim_{\varepsilon\to 0}\left(\lim_{N\to \infty} \frac{\log \#(X, d_N, \varepsilon)}{N}\right), \\  
   \mdim(X, T) & = \lim_{\varepsilon\to 0}\left(\lim_{N\to \infty} \frac{\widim_\varepsilon(X, d_N)}{N}\right).
\end{align*}
These are topological invariants, i.e., their values are independent of the choice of $d$.
Mean dimension was introduced by Gromov \cite{Gromov99}.
We also define the \textbf{lower metric mean dimension} $\lmdimm(X,T, d)$ by 
\[ \lmdimm(X, T, d) =
     \liminf_{\varepsilon\to 0} \left(\lim_{N\to \infty} \frac{\log  \#(X, d_N, \varepsilon)}{N\log (1/\varepsilon)}\right). \]
This is a metric-dependent quantity.
This notion was introduced by Lindenstrauss--Weiss \cite{Lindenstrauss--Weiss, Lindenstrauss}.
They proved that $\mdim(X, T) \leq \lmdimm(X, T, d)$ \cite[Theorem 4.2]{Lindenstrauss--Weiss}.
In particular, if $\htop(X, T)$ is finite then $\mdim(X, T)$ is zero.

Let $(X, T)$ and $(Y, S)$ be dynamical systems.
Let $\pi\colon X\to Y$ be an equivariant continuous map. 
Here “equivariant” means that $\pi\circ T = S\circ \pi$.
We often denote it by $\pi\colon (X, T) \to (Y, S)$ for clarifying the underlying dynamics.
Let $d$ be a metric on $X$.
We define the \textbf{topological conditional entropy} of $\pi$ by 
\[  \htop(\pi, T) 
    = \lim_{\varepsilon \to 0} \left(\lim_{N\to \infty} \frac{\sup_{y\in Y} \log \#\left(\pi^{-1}(y), d_N, \varepsilon\right)}{N}\right). \]
This is also a topological invariant.
This quantity measures how much entropy exists in the fibers of $\pi$.
See the book of Downarowicz \cite[Section 6]{Downarowicz} for the details about topological conditional entropy.
We also define the \textbf{lower conditional metric mean dimension} of $\pi$ by 
\[  \lmdimm(\pi, T, d) =   \liminf_{\varepsilon\to 0}  
     \left(\lim_{N\to \infty} \frac{\sup_{y\in Y} \log \#\left(\pi^{-1}(y), d_N, \varepsilon\right)}{N\log (1/\varepsilon)}\right). \]
This notion was introduced by Liang \cite{Liang}.
If $\lmdimm(\pi, T, d)$ is positive then $\htop(\pi, T)$ is infinite.
Hence we can say that $\lmdimm(\pi, T, d)$ quantifies \textit{how infinite $\htop(\pi, T)$ is}.

For an equivariant continuous map $\pi\colon (X, T)\to (Y, S)$, Bowen \cite[Theorem 17]{Bowen71} proved that 
\begin{equation}  \label{eq: Bowen inequality}
     \htop(X, T) \leq  \htop(Y, S) + \htop(\pi, T).
\end{equation}
Intuitively, this inequality says that if $(X, T)$ is \textit{larger than} $(Y, S)$ then some fiber of $\pi$ must be \textit{large}.
A motivation of \cite{Shi--Tsukamoto} is to study a situation where Bowen’s inequality (\ref{eq: Bowen inequality}) becomes meaningless.
Namely, the paper \cite{Shi--Tsukamoto} studied the situations where both $\htop(X, T)$ and $\htop(Y, S)$ are infinite.
In this case the inequality (\ref{eq: Bowen inequality}) provides no information about the fibers of $\pi$.

Let $a$ be a natural number. We consider the bi-infinite product of the $a$-dimensional cube $[0, 1]^a$:
\[  \left([0, 1]^a\right)^{\mathbb{Z}} = \cdots \times [0,1]^a\times [0,1]^a\times [0,1]^a\times \cdots. \]
We define the shift map $\sigma\colon \left([0,1]^a\right)^{\mathbb{Z}} \to \left([0,1]^a\right)^{\mathbb{Z}}$ by 
\[  \sigma\left((x_n)_{n\in \mathbb{Z}}\right) = (x_{n+1})_{n\in \mathbb{Z}}, \quad (\text{each } x_n\in [0,1]^a). \]
The dynamical system $\left(\left([0,1]^a\right)^{\mathbb{Z}}, \sigma\right)$ is called the \textbf{full-shift on the alphabet $[0,1]^a$}.
It has infinite topological entropy.
We denote by $\norm{\cdot}_\infty$ the max-norm of $\mathbb{R}^a$;
$\norm{(u_1, \dots, u_a)}_\infty = \max_{1\leq i \leq a}  |u_i|$.
We define a metric $D$ on $\left([0, 1]^a\right)^{\mathbb{Z}}$ by 
\[ D\left((x_n)_{n\in \mathbb{Z}}, (y_n)_{n\in \mathbb{Z}}\right) =  \sum_{n\in \mathbb{Z}} 2^{-|n|} \norm{x_n-y_n}_\infty. \]
The metric mean dimension of $\left(\left([0,1]^a\right)^{\mathbb{Z}}, \sigma, D\right)$ is equal to $a$.
Its (topological) mean dimension is also equal to $a$.

The following theorem was proved in \cite[Theorem 1.5]{Shi--Tsukamoto}.

\begin{theorem}[\cite{Shi--Tsukamoto}] \label{theorem: Shi--Tsukamoto}
Let $a$ be a natural number and let $(Y, S)$ be a dynamical system.
For any equivariant continuous map $\pi\colon \left(\left([0,1]^a\right)^{\mathbb{Z}}, \sigma\right) \to (Y, S)$ we have 
\begin{equation} \label{eq: lower bound on conditional metric mean dimension}
     \lmdimm(\pi, \sigma, D) \geq  a - 2\mdim(Y, S). 
\end{equation}     
Here $\mdim(Y, S)$ is the mean dimension of $(Y, S)$.
\end{theorem}

The factor $2$ in the right-hand side of (\ref{eq: lower bound on conditional metric mean dimension}) looks unsatisfactory.
Therefore the paper \cite[Problem 1.6]{Shi--Tsukamoto} asked whether we can remove it or not.
The purpose of the present paper is to solve this question affirmatively.
Our main result is the following.

\begin{theorem}[Main theorem]  \label{main theorem}
Let $a$ be a natural number and let $(Y, S)$ be a dynamical system.
For any equivariant continuous map $\pi\colon \left(\left([0,1]^a\right)^{\mathbb{Z}}, \sigma\right) \to (Y, S)$ we have 
\begin{equation}  \label{eq: improved lower bound on conditional metric mean dimension}
     \lmdimm(\pi, \sigma, D) \geq  a - \mdim(Y, S). 
\end{equation}     
In particular, if $a>\mdim(Y, S)$ then any equivariant continuous map 
$\pi\colon \left(\left([0,1]^a\right)^{\mathbb{Z}}, \sigma\right) \to (Y, S)$
has infinite topological conditional entropy.
\end{theorem}

We can easily see that the estimate (\ref{eq: improved lower bound on conditional metric mean dimension})
is optimal \cite[Example 1.4]{Shi--Tsukamoto}:
Let $a\geq b$ be two natural numbers and let $f\colon  [0,1]^a\to  [0,1]^b$ be the projection to the first 
$b$ coordinates.
Define $\pi\colon \left([0, 1]^a\right)^{\mathbb{Z}} \to \left([0,1]^b\right)^{\mathbb{Z}}$ by 
$\pi\left((x_n)_{n\in \mathbb{Z}}\right) = \left(f(x_n)\right)_{n\in \mathbb{Z}}$.
This is a shift equivariant continuous map between two full-shifts.
We have 
\[  \lmdimm(\pi, \sigma, D) = a-b. \]

The following question was posed in \cite[Problem 1.7]{Shi--Tsukamoto}.
We would like to repeat it here.

\begin{problem}[\cite{Shi--Tsukamoto}]
Let $(X, T)$ and $(Y, S)$ be dynamical systems with $\mdim(X, T) > \mdim(Y, S)$.
Does every equivariant continuous map $\pi\colon (X, T)\to (Y, S)$ have infinite topological conditional entropy?
\end{problem}

Theorem \ref{main theorem} affirmatively solves this problem in the case that $(X, T)$ is the full-shift on the alphabet $[0,1]^a$.
But the problem is widely open in general.

The basic structure of the proof of Theorem \ref{main theorem} is the same as that of Theorem \ref{theorem: Shi--Tsukamoto}
given in \cite{Shi--Tsukamoto}.
The (only one) additional ingredient is a classical theorem of Hurewicz \cite[p.124]{Kuratowski} 
which we will review in \S \ref{section: preparations} below.
The author was ignorant about this theorem at the time of writing \cite{Shi--Tsukamoto}.

\vspace{0.3cm}

\textbf{Organization of the paper.}
We prepare the waist inequality and the theorem of Hurewicz in \S \ref{section: preparations}.
We also prepare a useful result on conditional metric mean dimension in \S \ref{section: preparations}.
We prove Theorem \ref{main theorem} in \S \ref{section: proof of the main theorem}.
We explain the proof of (a corollary of) the theorem of Hurewicz in the appendix for completeness.

\vspace{0.3cm} 

\textbf{Acknowledgement.}
I would like to thank Professor Michael Levin.
I learnt the theorem of Hurewicz (Theorem \ref{theorem: Hurewicz} below) from him.
Recently he proved a very interesting mean dimension version of Theorem \ref{theorem: Hurewicz} in \cite{Levin}.

\section{Preparations}  \label{section: preparations}

\subsection{Klartag’s waist inequality}  \label{subsection: Klartag’s waist inequality}

Waist inequality was first discovered by Gromov \cite{Gromov03}.
It has several variations. (See \cite{Guth}.)
Here we review a convex geometry version of waist inequality developed by Klartag \cite{Klartag17}.

For $r>0$ and subsets $A, B$ of $\mathbb{R}^n$ we set $A+B = \{x+y\mid x\in A, y\in B\}$
and $rA = \{rx\mid x\in A\}$.
A \textbf{convex body} of $\mathbb{R}^n$ is a compact convex subset of $\mathbb{R}^n$ with non-empty interior.
A convex body $K$ is said to be \textbf{centrally-symmetric} if $-K = K$ where $-K := \{-x\mid x\in K\}$.
A function $\varphi\colon \mathbb{R}^n \to [0, \infty)$ is said to be \textbf{log-concave} if for any $x, y\in \mathbb{R}^n$ and 
$0<\lambda < 1$ we have $\varphi\left(\lambda x + (1-\lambda) y\right) \geq  \varphi(x)^{\lambda}\varphi(y)^{1-\lambda}$.

Klartag \cite[Theorem 5.7]{Klartag17} proved the following version of the waist inequality.

\begin{theorem}[Klartag 2017] \label{theorem: Klartag}
Let $m\leq n$ be natural numbers, and let $K\subset \mathbb{R}^n$ be a centrally-symmetric convex body.
Let $\mu$ be a probability measure supported in $K$ which has a log-concave density function with respect to the 
Lebesgue measure. Then for any continuous map $f\colon K\to \mathbb{R}^m$ there exists $t\in \mathbb{R}^m$ satisfying 
\[   \mu\left(f^{-1}(t)+rK\right) \geq  \left(\frac{r}{2+r}\right)^m  \quad \text{for all $0<r<1$}. \]
\end{theorem}

We denote the max norm of $\mathbb{R}^n$ by $\norm{\cdot}_\infty$: 
$\norm{(u_1, \dots, u_n)}_\infty = \max_{1\leq k \leq n} |u_k|$.
For $r>0$ and a compact subset $A$ of $\mathbb{R}^n$ we define $\#\left(A, \norm{\cdot}_\infty, r\right)$ as the minimum number $N$
such that there exists an open cover $\{U_1, \dots, U_N\}$ of $A$ satisfying 
$\diam (U_i, \norm{\cdot}_\infty) < r$ for all $1\leq i \leq  N$.
If $A$ is empty, we set $\#\left(A, \norm{\cdot}_\infty, r\right) =0$.

\begin{corollary}  \label{cor: corollary of Klartag waist inequality}
Let $m\leq n$ be two natural numbers.
For any continuous map $f\colon [0, 1]^n\to \mathbb{R}^m$ there exists a point $t\in \mathbb{R}^m$
satisfying 
\[  \#\left(f^{-1}(t), \norm{\cdot}_\infty, r\right) \geq \frac{1}{8^n}\left(\frac{1}{r}\right)^{n-m} \quad 
     \text{for all } 0< r< \frac{1}{2}. \]
\end{corollary}

\begin{proof}
This is a simple consequence of Klartag’s waist inequality (Theorem \ref{theorem: Klartag}).
See \cite[Corollary 2.4]{Shi--Tsukamoto} for a proof.
\end{proof}

\subsection{A theorem of Hurewicz}    \label{subsection: a theorem of Hurewicz}

Let $(X, d)$ be a compact metric space.
We define the \textbf{topological dimension} of $X$ by 
$\dim X = \lim_{\varepsilon\to 0} \widim_\varepsilon(X, d)$.
This is a topological invariant.
Let $n$ be a natural number. We define $C(X, [0,1]^n)$ as the space of all continuous maps from $X$ to $[0,1]^n$.
This space has a metric defined by 
\[  \rho(f, g) = \sup_{x\in X} |f(x)-g(x)|, \quad (f, g\in C(X, [0,1]^n)). \]
$\left(C(X, [0,1]^n), \rho\right)$ is a complete metric space.
In particular, we can apply to it the Baire category theorem.

The following classical result of Hurewicz \cite[p.124]{Kuratowski} is an essential ingredient of the proof of Theorem \ref{main theorem}.

\begin{theorem} \label{theorem: Hurewicz}
Suppose $n>\dim X$. 
There exists a dense $G_\delta$ subset $A$ of $C(X, [0,1]^n)$ such that every $f\in A$ satisfies
\[   \forall y\in [0,1]^n:  \quad |f^{-1}(y)| \leq  \left\lceil \frac{\dim X + 1}{n -\dim X}\right\rceil. \]
Here $|f^{-1}(y)|$ is the cardinality of $f^{-1}(y)$.
\end{theorem}

In particular, letting $n = \dim X + 1$, we get the following corollary.

\begin{corollary} \label{cor: Hurewicz}
 Let $X$ be a finite dimensional compact metrizable space. 
There exists a dense $G_\delta$ subset $A$ of $C\left(X, [0,1]^{\dim X + 1}\right)$ such that every $f\in A$ satisfies
\[  \forall y \in [0,1]^{\dim X +1}: \quad |f^{-1}(y)| \leq \dim X + 1. \]  
\end{corollary}

For example, if $X$ is one dimensional then there is a continuous map $f\colon X\to [0,1]^2$ such that every fiber of $f$ contains at most 
two points. If $X$ is a graph, this is intuitively obvious.
Corollary \ref{cor: Hurewicz} shows that this holds for all dimensions.

We use Corollary \ref{cor: Hurewicz} in the proof of Theorem \ref{main theorem}.
For completeness, we explain a proof of Corollary \ref{cor: Hurewicz} in the appendix.

\subsection{A formula of conditional metric mean dimension} \label{subsection: a formula of conditional metric mean dimension}

Let $(X, T)$ and $(Y, S)$ be dynamical systems.
Let $d$ and $d^\prime$ be metrics on $X$ and $Y$ respectively.
For each natural number $N$, we define $d_N$ and $d^\prime_N$ by 
$d_N(x_1, x_2) = \max_{0\leq n <N} d(T^n x_1, T^n x_2)$ and 
$d^\prime_N(y_1, y_2) = \max_{0\leq n < N} d^\prime(S^n y_1, S^n y_2)$. 
For $y\in Y$ and $r>0$ we define 
\[  B_r(y, d^\prime_N) = \{x\in Y\mid d^\prime_N(x, y) \leq r\}. \]

\begin{lemma}  \label{lemma: formula of metric mean dimension}
For any equivariant continuous map $\pi\colon (X, T)\to (Y, S)$ we have
\[  \lmdimm(\pi, T, d) = \liminf_{\varepsilon\to 0} \left\{\lim_{\delta\to 0} \left(\lim_{N\to \infty}
     \frac{\sup_{y\in Y} \log \#\left(\pi^{-1}\left(B_\delta(y, d^\prime_N)\right), d_N, \varepsilon\right)}{N\log (1/\varepsilon)}\right)\right\}. \]
\end{lemma}

\begin{proof}
See \cite[Lemma 3.2]{Shi--Tsukamoto}.
\end{proof}

\section{Proof of Theorem \ref{main theorem}}  \label{section: proof of the main theorem}

The purpose of this section is to prove Theorem \ref{main theorem}.
For potential applications in a future, we prove a more general statement as follows.
This is a refinement of \cite[Theorem 4.1]{Shi--Tsukamoto}.

\begin{theorem} \label{theorem: technical theorem}
Let $s\geq 0$ be a real number and let $(X, T)$ be a dynamical system with a metric $d$ on $X$.
Suppose that there exist sequences of natural numbers $N_n, M_n$ and continuous maps 
$\psi_n\colon [0,1]^{M_n}\to X$ $(n\geq 1)$ satisfying 
  \begin{itemize}
   \item $N_n\to \infty$ and $M_n/N_n\to s$ as $n\to \infty$,
   \item $\norm{x-y}_\infty \leq  d_{N_n}\left(\psi_n(x), \psi_n(y)\right)$ for all $x, y\in [0,1]^{M_n}$.
  \end{itemize}
  Then for any dynamical systems $(Y, S)$ and any equivariant continuous maps $\pi\colon (X, T)\to (Y, S)$ we have
  \[ \lmdimm(\pi, T, d) \geq   s- \mdim(Y, S). \]
\end{theorem}

We note that, under the assumption of this theorem, we have $\mdim(X, T) \geq s$ \cite[Theorem 4.1 (1)]{Shi--Tsukamoto}.

\begin{proof}
The proof closely follows the argument of \cite[Theorem 4.1]{Shi--Tsukamoto}.
The difference is that we use a theorem of Hurewicz (Corollary \ref{cor: Hurewicz}).

If $\mdim(Y, S)$ is infinite then the statement is trivial. So we assume that $\mdim(Y, S)$ is finite.
We take a metric $d^\prime$ on $Y$. (Any metric will work.)
Let $0<\varepsilon<1/2$ and $\delta>0$ be arbitrary numbers.
From the definition of mean dimension,
there exists a natural number $N(\delta)$ such that for any natural number $N\geq N(\delta)$ there 
exists a finite simplicial complex $K_N$ and a $\delta$-embedding 
$\varphi_N\colon  (Y, d^\prime_N)\to K_N$ satisfying 
\begin{equation} \label{eq: dimension of K_N}
  \frac{\dim K_N}{N} < \mdim(Y, S) + \delta.
\end{equation}
By Corollary \ref{cor: Hurewicz} there is a continuous map $g_N\colon K_N \to [0,1]^{\dim K_N +1}$ such that for every point 
$p \in [0, 1]^{\dim K_N +1}$ we have 
\begin{equation*}  
   |g_N^{-1}(p)|  \leq   \dim K_N  + 1.  
\end{equation*}   

For each $n$ with $N_n\geq N(\delta)$ we set $\pi_n = \varphi_{N_n}\circ \pi \circ \psi_n\colon [0,1]^{M_n}\to  K_{N_n}$.
We consider the map 
$g_{N_n}\circ  \pi_n \colon [0,1]^{M_n}\to  [0,1]^{\dim K_{N_n}+1}$.
We apply to it the waist inequality (Corollary \ref{cor: corollary of Klartag waist inequality}):
There exists a point $p \in [0,1]^{\dim K_{N_n}+1}$ (depending on $n$) such that 
\[  \#\left(\left(g_{N_n}\circ \pi_n\right)^{-1}(p), \norm{\cdot}_\infty, \varepsilon\right) \geq  \frac{1}{8^{M_n}}
     \left(\frac{1}{\varepsilon}\right)^{M_n -\dim K_{N_n}-1}. \]
Since we know that $|g_{N_n}^{-1}(p)|  \leq   \dim K_{N_n}  + 1$, there exists a point $t \in g_{N_n}^{-1}(p)$ such that 
\begin{equation} \label{eq: choosing a point t by waist and Hurewicz}
    \#\left(\pi_n^{-1}(t), \norm{\cdot}_\infty, \varepsilon\right) \geq  \frac{1}{8^{M_n} (\dim K_{N_n}+1)}
     \left(\frac{1}{\varepsilon}\right)^{M_n -\dim K_{N_n}-1}. 
\end{equation}
Take $t^\prime \in Y$ (which also depends on $n$) satisfying $\varphi_{N_n}(t^\prime) = t$.

\begin{claim}  \label{claim: t and t^prime}
  We have $\psi_n\left(\pi_n^{-1}(t)\right) \subset  \pi^{-1}\left(B_\delta(t^\prime, d^\prime_{N_n})\right)$. 
\end{claim}
\begin{proof}
Let $x\in \pi_n^{-1}(t)$. Then
\[  \varphi_{N_n}\left(\pi\left(\psi_n(x)\right)\right) = \pi_n(x) = t = \varphi_{N_n}(t^\prime). \]
Recall that the map $\varphi_{N_n}\colon Y\to K_{N_n}$ is a $\delta$-embedding with respect to $d^\prime_{N_n}$.
Hence $d^\prime_{N_n}\left(\pi\left(\psi_n(x)\right), t^\prime\right) < \delta$. Thus 
\[  \pi\left(\psi_n(x)\right)  \in  B_\delta(t^\prime, d^\prime_{N_n}). \]
Namely $\psi_n\left(\pi_n^{-1}(t)\right) \subset \pi^{-1}\left(B_\delta(t^\prime, d^\prime_{N_n})\right)$. 
\end{proof}

Recall that we have assumed $\norm{x_1-x_2}_\infty \leq d_{N_n}\left(\psi_n(x_1), \psi_n(x_2)\right)$ for all $x_1, x_2\in [0,1]^{M_n}$.
Therefore 
\begin{align*}
 \#\left(\pi_n^{-1}(t), \norm{\cdot}_\infty, \varepsilon\right) & \leq 
 \#\left(\psi_n\left(\pi_n^{-1}(t)\right), d_{N_n}, \varepsilon\right) \\
 & \leq \#\left(\pi^{-1}\left(B_\delta(t^\prime, d^\prime_{N_n})\right), d_{N_n}, \varepsilon\right).
\end{align*}
In the second inequality we have used Claim \ref{claim: t and t^prime}.
It follows from (\ref{eq: choosing a point t by waist and Hurewicz}) that
\[  \sup_{y\in Y} \#\left(\pi^{-1}\left(B_\delta(y, d^\prime_{N_n})\right), d_{N_n}, \varepsilon\right)  \geq 
     \frac{1}{8^{M_n} (\dim K_{N_n}+1)}
     \left(\frac{1}{\varepsilon}\right)^{M_n -\dim K_{N_n}-1}. \]
We take the logarithm of the both sides and divide them by $N_n \log (1/\varepsilon)$.
\begin{align*}
   &\frac{\sup_{y\in Y} \log  \#\left(\pi^{-1}\left(B_\delta(y, d^\prime_{N_n})\right), d_{N_n}, \varepsilon\right)}{N_n \log (1/\varepsilon)} \\
  & \geq  \frac{M_n-\dim K_{N_n}-1}{N_n} - \frac{M_n\log 8 + \log \left(\dim K_{N_n} +1\right)}{N_n \log (1/\varepsilon)} \\
  & \geq \frac{M_n}{N_n} - \mdim(Y, S) -\delta - \frac{1}{N_n} 
     - \frac{M_n\log 8 + \log \left\{1+ N_n(\mdim(Y, S) + \delta)\right\}}{N_n \log (1/\varepsilon)}.
\end{align*}  
In the second inequality we have used (\ref{eq: dimension of K_N}).
Recall $N_n\to \infty$ and $M_n/N_n\to s$ as $n\to \infty$. We let $n\to \infty$ and get 
\[ \lim_{N\to \infty}
 \frac{\sup_{y\in Y} \log  \#\left(\pi^{-1}\left(B_\delta(y, d^\prime_{N})\right), d_{N}, \varepsilon\right)}{N \log (1/\varepsilon)} 
 \geq  s - \mdim(Y, S)-\delta - \frac{s\log 8}{\log (1/\varepsilon)}. \]
Letting $\delta \to 0$ and $\varepsilon \to 0$, we conclude 
\[  \lim_{\varepsilon \to 0} \left\{\lim_{\delta \to 0}
     \left(\lim_{N\to \infty} \frac{\sup_{y\in Y} 
     \log  \#\left(\pi^{-1}\left(B_\delta(y, d^\prime_{N})\right), d_{N}, \varepsilon\right)}{N \log (1/\varepsilon)}\right)\right\}
     \geq  s - \mdim(Y, S). \]
The left-hand side is equal to $\lmdimm(\pi, T, d)$ by Lemma \ref{lemma: formula of metric mean dimension}.
Thus we have finished the proof.
\end{proof}

Now we can prove Theorem \ref{main theorem}.
Recall that the full-shift $\left([0,1]^a\right)^{\mathbb{Z}}$ has the metric 
\[  D\left((x_n)_{n\in \mathbb{Z}}, (y_n)_{n\in \mathbb{Z}}\right) = \sum_{n\in \mathbb{Z}}2^{-|n|}\norm{x_n-y_n}_\infty. \]

\begin{proof}[Proof of Theorem \ref{main theorem}]
Let $\pi$ be any equivariant continuous map from $\left(\left([0,1]^a\right)^{\mathbb{Z}}, \sigma\right)$ to 
a dynamical system $(Y, S)$. We want to show $\lmdimm(\pi, \sigma, D) \geq  a - \mdim(Y, S)$.
For natural numbers $n$ we define continuous maps $\psi_n\colon \left([0,1]^a\right)^n \to \left([0,1]^a\right)^{\mathbb{Z}}$ by
\[ \psi_n(x_0, x_1, \dots, x_{n-1})_k = \begin{cases} x_k \quad & (0\leq k \leq n-1) \\  0 \quad &(\text{otherwise}) \end{cases}. \]
We have $D_n\left(\psi_n(x), \psi_n(y)\right) \geq  \norm{x-y}_\infty$ for all $x, y\in \left([0,1]^a\right)^n = [0,1]^{an}$.
Thus we can apply Theorem \ref{theorem: technical theorem} to the present setting
with the parameters $N_n = n$, $M_n= an$ and $s = a$.
Then we conclude that $\lmdimm(\pi, \sigma, D) \geq  a - \mdim(Y, S)$.
\end{proof}

\appendix

\section{Proof of Corollary \ref{cor: Hurewicz}}

We explain a proof of Corollary \ref{cor: Hurewicz} here for convenience of readers.
The proof of Theorem \ref{theorem: Hurewicz} (the theorem of Hurewicz) is similar 
but we concentrate on the statement of Corollary \ref{cor: Hurewicz}
for simplicity.
Our proof follows \cite[p.124]{Kuratowski}.

Let $P$ be a simplicial complex. 
We denote by $V(P)$ the set of vertices of $P$. 
We also denote by $V(\Delta)$ the set of vertices of $\Delta$ for each simplex $\Delta \subset P$.
(Here “simplex” means a closed simplex.)
A map $g\colon P\to \mathbb{R}^m$ is called a \textbf{simplicial map} if for every simplex $\Delta\subset P$
of vertices $v_0, v_1, \dots, v_n$ we have 
\[ g\left(\sum_{i=0}^n t_i v_i\right) = \sum_{i=0}^n t_i g(v_i), \]
where $t_i\geq 0$ and $t_0+t_1+\dots+t_n=1$.
A simplicial map $g$ is determined by the values $g(v)$ $(v\in V(P))$.

The next lemma was given in \cite[Lemma 2.1]{Gutman--Lindenstrauss--Tsukamoto}

\begin{lemma} \label{lemma: simplicial approximation}
Let $(X, d)$ be a compact metric space and let $f\colon X\to \mathbb{R}^m$ be a continuous map.
Let $\delta$ and $\varepsilon$ be positive numbers such that 
\[  d(x, y) < \varepsilon  \Longrightarrow  |f(x)-f(y)| < \delta. \]
Let $\pi\colon X\to P$ be an $\varepsilon$-embedding to some finite simplicial complex $P$.
Then, after subdividing $P$ sufficiently finer, we can find a simplicial map $g\colon P\to \mathbb{R}^m$
satisfying $|f(x)-g(\pi(x))| < \delta$ for all $x\in X$.
\end{lemma}

\begin{proof}
For each vertex $v\in V(P)$ we define the \textbf{open star} $\mathrm{Star}(v)$ as the union of the interior of 
all simplices of $P$ containing $v$.
(For a simplex $\Delta\subset P$ of vertices $v_0, \dots, v_k$, 
we define its \textbf{interior} as the set of $\sum_{i=0}^k \lambda_i v_i$
where $\lambda_i >0$ for all $i$ and $\sum_i \lambda_i  = 1$. In particular, the interior of $\{v\}$ is $\{v\}$.)

Since $\pi$ is an $\varepsilon$-embedding, we can subdivide $P$ sufficiently fine so that we have 
$\diam \pi^{-1}\left(\mathrm{Star}(v)\right) < \varepsilon$ for all $v\in V(P)$.

Let $v\in V(P)$. If $\pi^{-1}\left(\mathrm{Star}(v)\right) =\emptyset$ then we set $g(v) = 0$.
If $\pi^{-1}\left(\mathrm{Star}(v)\right) \neq \emptyset$ then we pick up $x_v\in \pi^{-1}\left(\mathrm{Star}(v)\right)$ and 
set $g(v) = f(x_v)$. We extend $g$ over $P$ by linearity.

Let $x\in X$. We choose a simplex $\Delta \subset P$ whose interior contains $\pi(x)$.
We have $x\in \pi^{-1}\left(\mathrm{Star}(v)\right)$ for all $v\in V(\Delta)$ and hence 
$|f(x) - f(x_v)| < \delta$.
Therefore we have $|f(x) - g(\pi(x))| < \delta$.
\end{proof}

\begin{lemma} \label{lemma: linear algebra}
Let $\delta>0$ and let $P$ be a finite simplicial complex of dimension $n$.
For any simplicial map $f\colon P\to \mathbb{R}^{n+1}$ there exists a simplicial map 
$g\colon P\to \mathbb{R}^{n+1}$ satisfying 
   \begin{itemize}
      \item  $|f(x)-g(x)| < \delta$ for all $x\in P$,
      \item  for any pairwisely disjoint simplices $\Delta_1, \dots, \Delta_{n+2}$ of $P$ we have 
                $g(\Delta_1)\cap \dots \cap g(\Delta_{n+2}) = \emptyset$.    
   \end{itemize}
\end{lemma}

\begin{proof}
Let $W_1$ and $W_2$ be affine subspaces of $\mathbb{R}^{n+1}$.
We say that $W_1$ and $W_2$ are \textbf{in general position} if the following two conditions\footnote{The condition (2) looks complicated.
We can understand it in terms of the projective space. Let $\mathbb{P}^{n+1}$ be the real projective space of dimension $n+1$.
We naturally consider $\mathbb{R}^{n+1}$ as a subset of $\mathbb{P}^{n+1}$.
For an affine subspace $W$ of $\mathbb{R}^{n+1}$ we denote by $\overline{W}$ the closure of $W$ in $\mathbb{P}^{n+1}$.
Then the condition (2) means that if $\dim W_1 + \dim W_2 < n+1$ then $\overline{W_1}$ and $\overline{W_2}$ do not intersect in $\mathbb{P}^{n+1}$.}
are satisfied:
\begin{enumerate}
   \item If $\dim W_1 + \dim W_2\geq n+1$ then $\dim (W_1\cap W_2) = \dim W_1 + \dim W_2 -n-1$.
   \item If $\dim W_1 + \dim W_2 < n+1$ then $W_1\cap W_2 = \emptyset$ and, moreover, for all $u_1 \in W_1$ and $u_2\in W_2$
    we have $(-u_1+W_1) \cap (-u_2+W_2) = \{0\}$.
\end{enumerate} 
The latter condition of (2) means that any line on $W_1$ is not parallel to any line on $W_2$.
(Also, notice that if we have $(-u_1+W_1) \cap (-u_2+W_2) = \{0\}$ for \textit{some} $u_1\in W_1$ and $u_2\in W_2$ then 
we have the same condition for \textit{all} $u_1\in W_1$ and $u_2\in W_2$.)
If $W_1$ and $W_2$ are in general position then their small perturbations are also in general position.
Given any affine subspaces $W_1$ and $W_2$ of $\mathbb{R}^{n+1}$, we can perturb $W_1$ arbitrarily small so that 
$W_1$ and $W_2$ are in general position.

For each $v\in V(P)$ we can choose $g(v)\in \mathbb{R}^{n+1}$ so that 
\begin{itemize}
  \item $|g(v)-f(v)| < \delta$ for all $v\in V(P)$,
  \item $g|_\Delta\colon \Delta \to \mathbb{R}^{n+1}$ is injective for every simplex $\Delta \subset P$ where 
  we define the map $g$ on $\Delta$ by linearity.
\end{itemize}
Notice that these conditions are stable under small perturbations of $g(v)$.

For a simplex $\Delta\subset P$ 
we set $W\left(g(\Delta)\right) = \{\sum_{v\in V(\Delta)} \lambda_v g(v)\mid  \sum_{v\in V(\Delta)} \lambda_v =1\}$.
This is the affine subspace of $\mathbb{R}^{n+1}$ spanned by $g(\Delta)$. We have 
$\dim W\left(g(\Delta)\right) = \dim g(\Delta) = \dim \Delta$.

By the induction on the number $m$ of simplices, we can slightly perturb $g(v)$ $(v\in V(P))$ so that, for any pairwisely disjoint 
simplices $\Delta_1, \dots, \Delta_m$ of $P$, the subspaces $W\left(g(\Delta_m)\right)$ and 
$\bigcap_{i=1}^{m-1} W\left(g(\Delta_i)\right)$ are in general position.

In particular (recalling $\dim P = n$), we have $W\left(g(\Delta_1)\right) \cap \dots \cap W\left(g(\Delta_m)\right) = \emptyset$
for $m \geq n+2$.
Thus we have $g(\Delta_1) \cap \dots \cap g(\Delta_{n+2}) = \emptyset$ for any pairwisely disjoint simplices 
$\Delta_1, \dots, \Delta_{n+2}$ of $P$.
\end{proof}

We write the statement of Corollary \ref{cor: Hurewicz} again.

\begin{proposition}[$=$ Corollary \ref{cor: Hurewicz}]
Let $(X, d)$ be a finite dimensional compact metric space. Let $n= \dim X$.
There exists a dense $G_\delta$ subset $A$ of $C(X, \mathbb{R}^{n+1})$ such that every $f\in A$ satisfies
  \begin{equation} \label{eq: fiber cardinality}
     \forall y\in \mathbb{R}^{n+1}: \quad |f^{-1}(y)| \leq n+1. 
  \end{equation}   
\end{proposition}

\begin{proof}
Let $\mathcal{B}$ be a countable base of the topology of $X$.
Define 
\[ S = \{(U_1, \dots, U_{n+2})\in \mathcal{B}^{n+2}\mid \overline{U_i}\cap \overline{U_j} = \emptyset \> (i\neq j)\}. \]
This is a countable set. For $\mathbf{U} = (U_1, \dots, U_{n+2})\in S$ we set 
\[  C(\mathbf{U}) = \{f\in C(X, \mathbb{R}^{n+1}) \mid
      f\left(\overline{U_1}\right) \cap\dots \cap f\left(\overline{U_{n+2}}\right) = \emptyset\}. \]
We can easily see that this is an open subset of $C(X, \mathbb{R}^{n+1})$.
We will prove that it is also dense.
If this is proved, then the intersection $A:=\bigcap_{\mathbf{U}\in S} C(\mathbf{U})$ is a dense $G_\delta$ subset of $C(X, \mathbb{R}^{n+1})$
(by the Baire category theorem) and every $f\in \bigcap_{\mathbf{U}\in S} C(\mathbf{U})$ satisfies (\ref{eq: fiber cardinality}).

Let $\mathbf{U} = (U_1, \dots, U_{n+2})\in S$.
Let $f\in C(X, \mathbb{R}^{n+1})$ be an arbitrary map and let $\delta>0$.
We want to construct $f^\prime \in C(\mathbf{U})$ satisfying $|f(x)-f^\prime(x)|<\delta$ for all $x\in X$.

Take $\varepsilon>0$ such that 
\begin{itemize}
   \item  $d\left(\overline{U_i}, \overline{U_j}\right) > \varepsilon$ for $i\neq j$,
   \item  $d(x,y)< \varepsilon \Longrightarrow  |f(x)-f(y)| < \delta$.
\end{itemize}

From the definition of topological dimension, there exists an $\varepsilon$-embedding 
$\pi\colon X\to P$ to some finite simplicial complex $P$ with $\dim P \leq n$.
We have $\pi\left(\overline{U_i}\right) \cap \pi\left(\overline{U_j}\right) = \emptyset$ for $i\neq  j$.
Then we can subdivide $P$ sufficiently fine so that if two simplices $\Delta$ and $\Delta^\prime$ of $P$ satisfy 
$\Delta \cap \pi\left(\overline{U_i}\right) \neq \emptyset$ and $\Delta^\prime\cap \pi\left(\overline{U_j}\right) \neq \emptyset$
for some $i\neq j$ then $\Delta \cap \Delta^\prime =\emptyset$.

By Lemma \ref{lemma: simplicial approximation}, after subdividing $P$ further, there is a simplicial map 
$g\colon P\to \mathbb{R}^{n+1}$ satisfying $|f(x)-g(\pi(x))|< \delta$ for all $x\in X$.

By Lemma \ref{lemma: linear algebra}, we can perturb the map $g$ arbitrarily small so that for any pairwisely disjoint simplices 
$\Delta_1, \dots, \Delta_{n+2}$ of $P$ we have 
$g(\Delta_1) \cap \dots \cap g(\Delta_{n+2}) = \emptyset$.

Then the map $f^\prime:=g\circ \pi\colon X\to \mathbb{R}^{n+1}$ satisfies $|f(x)-f^\prime(x)|<\delta$ for all $x\in X$ and 
$f^\prime\left(\overline{U_1}\right) \cap\dots \cap f^\prime\left(\overline{U_{n+2}}\right) = \emptyset$.
\end{proof}

\end{document}